\documentclass[11pt]{article}
\usepackage[a4paper,margin=1in]{geometry}
\usepackage[T1]{fontenc}\usepackage[utf8]{inputenc}
\usepackage{lmodern,amsmath,amssymb,booktabs,array,microtype}
\usepackage[hidelinks]{hyperref}
\title{\textbf{Endpoint Architectures for Navier--Stokes}}
\author{Rishad Shahmurov\\\small Cellular Products Research and Development, Roswell, GA, USA}
\date{September 2026}
\begin{document}\maketitle
\begin{abstract}
We compare two endpoint analyses for three-dimensional Navier--Stokes flow: a general suitable-weak framework and the axisymmetric equations with swirl.  The common backend is selection of physical space--time records, nonoverlapping allocation of positive quantities, retention of the full localized state, and passage either to a compact transition law or to an explicitly named loss of compactness.  The front-end theorems are different.  In the general problem, bounded-critical sequences lead to symmetry/defect alternatives, while amplitude-fast normalization leads conditionally to an ancient Euler profile only after frequency tightness and a nonvanishing witness are retained.  In the axisymmetric problem, logarithmic swirl records are routed to meridional energy or scale-critical signed compression, whereas first-use provenance and recharge packing force recurrent fresh service toward boundary, age, frequency, or other recession mechanisms.  We state the interfaces and the remaining endpoints without identifying their currencies or reducing the general problem to the axisymmetric class.  No unconditional regularity theorem or singularity construction is claimed.
\end{abstract}

\section{Introduction}
The global regularity problem for the three-dimensional incompressible Navier--Stokes equations asks whether a smooth finite-energy solution can lose regularity in finite time.  Leray's construction of global weak solutions placed compactness at the center of the subject, while the partial-regularity theory of Caffarelli, Kohn, and Nirenberg showed that any singular set of a suitable weak solution is strongly constrained \cite{Leray1934,CKN1982}.  Geometric criteria such as the vorticity-direction condition of Constantin and Fefferman show that the organization of the flow, and not only its size, matters near a possible singularity \cite{ConstantinFefferman1993}.  Blow-up analysis therefore has to answer two questions simultaneously: what survives under rescaling, and whether the surviving object still carries the feature that made the original sequence singular.

A useful endpoint procedure begins with a physically selected space--time event rather than an arbitrary weakly convergent subsequence.  The selected event carries velocity and one fixed pressure representative, but also nonlinear products, signed fluxes, boundary or tail information, and the marks used to identify the event.  A positive quantity that has already been charged to one selected event should not be silently counted again at a later overlapping event.  The threshold decomposition of vorticity stretching provides one exact mechanism for such nonoverlapping assignment \cite{ShahmurovThreshold2026}.  Once a complete state has been selected, one faces a compactness alternative: either a subsequence converges in the topology needed by the equation, or a specific quantity escapes through scale, frequency, space, time, pressure, boundary flux, or amplitude.

The present paper compares two companion implementations of this philosophy.  The first is genuinely three-dimensional and works with suitable weak solutions and Euclidean symmetry defects \cite{ShahmurovEndpoint2026}.  The second is restricted to axisymmetric flows with swirl.  Its critical-face analysis uses the variables
\[
F=\frac{u^\theta}{r},\qquad
G=\frac{\omega^\theta}{r},\qquad
\Gamma=ru^\theta,
\]
and the exact source--Hodge structure of the axisymmetric equations \cite{ShahmurovAxisym2026}.  A separate companion paper develops first-use provenance, recharge packing, and the recession alternatives that arise when one asks for repeated fresh service from the same unforced solution \cite{ShahmurovAxisymProvenance2026}.

The distinction between these two axisymmetric papers is important.  There is a \emph{record-profile} problem and a \emph{recurrence/provenance} problem.  A large Type-II kinetic record may be meridional, swirl dominated, or tied to a critical signed-compression action.  That classification is not the same as asking whether an old current can legitimately perform a future service without prior use, and neither question can be renamed into the other.  The provenance analysis begins only after the critical-face and finite-depth preparation mechanisms have been separated.

The general three-dimensional analysis also has two scales.  On bounded-critical sequences one can seek a compact terminal state and test retained defect, symmetry, pressure, and nonlinear-product information.  On high-critical sequences one normalizes by a large kinetic amplitude.  The effective viscosity then tends to zero and the backward normalized time interval grows.  Under localized frequency tightness, pressure control, and a fixed nonvanishing witness, compactness is strong enough to pass the quadratic term and produce a nonzero ancient Euler profile.  Those hypotheses are not automatic consequences of the kinetic record: frequency/profile concentration and loss of the witness remain genuine alternatives.  Nor does an Euler limit by itself retain a singularity.

The axisymmetric structure is more rigid but not closed.  The circulation satisfies a homogeneous drift--diffusion equation and a maximum principle, the meridional velocity is recovered through a five-dimensional Hodge representation, and the source of $G$ is $\partial_z(F^2)$.  These identities yield sharp record inequalities and finite-history provenance formulas unavailable in a general three-dimensional flow.  They also show why compact, repeatedly first-use recurrence is difficult: represented fresh stock is depleted after use, while new stock must be assigned to ancestry, source generation, turnover, or boundary entry.  The analysis classifies the remaining recession mechanisms but does not eliminate all of them.

The common lesson is therefore methodological rather than identificatory.  Both programs need a complete same-solution state, a physically correct quotient, nonoverlapping ownership, and a closed transition relation.  But the endpoint theorem is different in each setting.  General three-dimensional flow must retain a singular or deformation witness through compactness and possible Euler passage.  Axisymmetric record analysis must route large swirl energy to a receiving branch.  Axisymmetric recurrence analysis must decide whether a genuinely new owner can enter from a noncompact boundary or age/frequency face and service another event.  None of these implications follows from the other two.

The organization reflects this separation.  Section 2 records the common threshold/transition backend.  Section 3 summarizes the general three-dimensional bounded- and high-critical endpoints.  Sections 4 and 5 separate the axisymmetric record and provenance fronts.  Section 6 compares the common structure and the nontransferable variables.  Section 7 states the current endpoints and limitations.

\section{Threshold selection and closed transition laws}
The threshold decomposition in \cite{ShahmurovThreshold2026} resolves finite-time vorticity stretching according to the amplitude level $a=|\omega|$.  For almost every threshold the stretching above $a$ equals an endpoint increment plus two nonnegative viscous terms: directional dissipation and level-surface amplitude dissipation.  Integrating in the threshold recovers the classical $L^p$ vorticity balance; at $p=2$ the two pieces reconstruct spacetime palinstrophy, while at $p=3/2$ they can be written through the scale-critical nonlinear root $|\omega|^{-1/4}\omega$.

The same positive threshold measure gives a first-assignment rule for overlapping selected events.  Under compactness, the normalized states need not be generated by one continuous deterministic return map.  It is enough that consecutive states remain in a compact chamber and belong to a closed physical transition relation.  Empirical measures of consecutive pairs then have subsequential limits with equal first and second marginals; disintegration yields a stationary Markov transition law supported on the closed relation.  This compactification principle neither creates compactness nor implies an atomic or self-similar recurrent profile.

This backend is intentionally weak enough to be shared by the two endpoint programs.  It provides accounting and a language for recurrence, but it does not identify the correct state variables for either the general or the axisymmetric equations.

\section{Three-dimensional endpoint structure}
The three-dimensional companion paper keeps a suitable-weak velocity--pressure state together with first-assigned nonnegative quantities, signed nonlinear currents, pressure normalization, tails, time faces, and Euclidean Killing defects \cite{ShahmurovEndpoint2026}.  The purpose of this enlarged state is to prevent a limiting profile from being declared regular or symmetric after the information selecting the singular sequence has disappeared.

On bounded-critical sequences, terminal-class quantization leads to a structured alternative.  A positive retained defect may survive; vanishing of the full Killing defect can, under the stated receiving interface, lead to an axisymmetric capture; or one is left with a nonsmooth terminal state separated by a positive Killing-defect floor.  The last branch cannot be removed from the size of the rotational defect alone.  The companion paper gives a divergence-free angular field whose rotationally averaged quadratic stress is isotropic, so a positive non-axisymmetric $L^2$ distance does not imply a quantitatively productive Reynolds forcing after projection.  Singularity-specific dynamics are still required on the positive-Killing-gap branch.

The high-critical branch uses a different normalization.  If $A_j$ is the selected kinetic amplitude, the amplitude-fast sequence has effective viscosity tending to zero and an increasingly long backward interval.  Under a local $L^2$ envelope, localized high-frequency tightness, pressure control, a negative-Sobolev time-derivative bound, and one fixed resolved spacetime witness, compactness yields strong local $L^2$ convergence.  The quadratic term then passes in $L^1_{\rm loc}$ and the limit is a nonzero ancient Euler solution.  Schematically,
\[
\boxed{
\begin{gathered}
\text{frequency or profile concentration, or witness loss}\
\vee\ \,\text{nonzero marked ancient Euler entrance}
\end{gathered}}
\]
The theorem is conditional in exactly the displayed sense: the singular record does not yet furnish the frequency-tightness or witness-retention hypotheses automatically.

A further spectral analysis concerns carriers that are null along an entire heat orbit.  Exact heat-orbit nullity separates on radial interaction shells; for real fields the shellwise covariance is isotropic, and normalized heat aging sends the $L^2$ mass out of every fixed spatial ball.  These facts constrain compact spectral carriers.  Their use in an actual singularity sequence remains downstream of the missing same-solution weld that produces the heat-orbit relation and retains the original singular/deformation witness.

Thus the high-critical endpoint is not ``Euler'' in a terminal sense.  If a classical Euler background is encountered while the singular witness survives, a further sub-viscosity normalization can remove that smooth background and force a deeper Navier--Stokes fluctuation.  Without retained singular information, however, the same statement is invalid.  The live alternatives remain a nonsmooth/noncompact Euler entrance, an infinite scale-orthogonal Navier--Stokes relay, or one of the explicit concentration/witness-loss exits.

\section{Axisymmetric Type-II records}
For axisymmetric flow with swirl, the critical-structure paper exploits the distinct scaling degrees of $F$, $G$, and $\Gamma$ \cite{ShahmurovAxisym2026}.  Besides mesoscopic and temporal critical faces, source/contact estimates, finite-depth preparation obstructions, and compactness results, the current record analysis gives two useful routes for large kinetic records.

Let
\[
A_R^2=\frac1R\int_{B_R}|u|^2,\qquad
M_R=\frac1R\int_{B_R}|u^{\rm mer}|^2,\qquad
S_R=A_R^2-M_R,
\]
for an axis-centered ball, and assume $|\Gamma|\le G$.  A sharp radial estimate gives
\[
S_R\le2\pi G^2\log\!\left(1+\frac{R D_{\theta,R}}{2\pi G^2}\right),
\]
where $D_{\theta,R}$ is controlled by the physical dissipation.  Consequently
\[
\int_I\left(e^{S_R/(2\pi G^2)}-1\right)\,dt
\le \frac{E_0R}{2\pi\nu G^2}
\]
for every time interval $I$.  On a parabolic interval of length $\tau R^2$, Jensen's inequality yields
\[
\langle S_R\rangle
\le2\pi G^2\log\!\left(1+\frac{E_0}{2\pi\nu G^2\tau R}\right).
\]
Thus a parabolic-time averaged total kinetic record with leading coefficient $C G^2\log(R_0/R)$ must contain a logarithmic meridional contribution whenever $C>2\pi$; in the quarter-meridional convention the meridional fraction exceeds one quarter for $C>8\pi/3$.  The same occupation estimate shows that, in any positive-density set of logarithmically high records above this threshold, the swirl-dominated subpopulation has vanishing parabolic density.  These statements require a genuine time-density hypothesis and do not turn isolated record times into a parabolic essential-supremum bound.

There is also an eventwise Type-II fork.  If a selected kinetic record is not already at least quarter-meridional, a cone estimate and the exact $q=2$ critical swirl balance force the scale-critical signed-compression action
\[
\mathcal A_2(t_0,t)
\ge
\frac{A_R^2}{24\pi\Gamma_*^2}
+\frac23\log A_R
-\frac13\log Z_2(t_0)-\log2.
\]
Hence a large record cannot hide indefinitely in an exponentially thin swirl cone while all critical compression output remains silent.  The clean structural fork is
\[
\boxed{
\text{meridional Type-II record}
\quad\vee\quad
\text{scale-critical signed-compression output}.}
\]
Neither receiving branch is excluded by this statement.  In particular, the record theorem is not an unconditional regularity criterion.

\section{Axisymmetric first-use provenance and recession}
The recurrence problem is developed separately in \cite{ShahmurovAxisymProvenance2026}.  It begins with the observation that a future service cannot be described only by the instantaneous values of $F$ and $G$: one must retain where the serving current came from and whether it has already been used.  The actual-drift backward trace admits an exact conjugation, and a common-cut extraction retains age distribution, regular-Lagrangian-flow provenance, ordered $P/Q$ data, cargo, and the unchanged full axisymmetric solution germ.

At the Kelvin-response level, the bookkeeping has four physical payers: ancestry already present at the cut, the native $\Gamma^2$ source, relative turnover, and boundary flux.  Represented fresh stock contracts after use, while genuinely new fixed-margin first-use blocks have a positive packing cost.  Therefore a finite represented and first-owned inventory cannot support infinitely many fixed-margin fresh recharges without leaving the compact represented regime.

The next step is a recession classification rather than a global contradiction.  Atomic service either second-normalizes to a qualifying ancient profile or fires a typed defect.  Diffuse service cannot remain a silent interior covariance and is routed toward source/cargo concentration, rate or coefficient recession, lower-face or boundary entrance, or tail/pressure-time loss.  The conclusion is deliberately one-sided:
\[
\boxed{
\begin{gathered}
\text{compact/passive fresh recurrence is exhausted}\\
\Longrightarrow\\
\text{tight ancient/classical exit or typed noncompact recession}
\end{gathered}}
\]
The typed alternatives are classified, not all eliminated.  In particular, boundary-fed or age/frequency-noncompact changing-owner service remains a genuine endpoint.  This is the current first-use problem and it should not be conflated with the Type-II record fork of the preceding section.

\section{Common structure and essential differences}
The two endpoint analyses share a backend that can be summarized as
\[
\begin{gathered}
\text{physical selection}\\
\longrightarrow\text{nonoverlapping ownership}\\
\longrightarrow\text{complete-state compactness or typed recession}\\
\longrightarrow\text{closed same-solution transition law}.
\end{gathered}
\]
The meaning of each state coordinate is nevertheless equation dependent.

\begin{center}
\small
\begin{tabular}{@{}p{0.16\textwidth}p{0.34\textwidth}p{0.34\textwidth}@{}}
\toprule
& General three-dimensional endpoint & Axisymmetric endpoint with swirl\\
\midrule
Selection & suitable-weak kinetic and defect records with fixed pressure gauge & axis-centered critical records and first-use service events\\[0.3em]
Main geometry & Euclidean Killing defect, pressure and product passage, singular/deformation witness & circulation, five-dimensional Hodge recovery, source $\partial_z(F^2)$, age and provenance\\[0.3em]
Compact terminal & retained defect; axisymmetric capture; nonsmooth positive-Killing-gap state & tight ancient/classical exit after finite-depth and recharge exhaustion\\[0.3em]
Noncompact terminal & frequency or profile concentration; witness loss; nonsmooth Euler entrance; infinite NS relay & boundary, age, frequency, source, cargo, tail, or pressure-time recession\\[0.3em]
Independent record front & amplitude-fast kinetic normalization & meridional Type-II record or critical signed compression\\
\bottomrule
\end{tabular}
\end{center}

Several tempting identifications are therefore invalid.  A general three-dimensional Killing defect is not axisymmetric circulation provenance.  An axisymmetric signed-compression action is not a general three-dimensional currency.  A stationary transition law need not contain an atom or a deterministic profile.  A nonzero ancient Euler entrance need not preserve the singular witness that selected the Navier--Stokes sequence.  Finally, the axisymmetric maximum principle and five-dimensional Hodge identities have no direct analogue that would reduce the general problem to the axisymmetric class.

\section{Current endpoints and limitations}
The synthesis now terminates at four explicit open fronts.  On bounded-critical general three-dimensional sequences, the nonsmooth positive-Killing-gap terminal remains.  On the high-critical side, one must either prove the frequency/witness compactness needed for a nonzero marked Euler entrance and then retain the singular mechanism, or classify the surviving nonsmooth Euler/infinite Navier--Stokes relay.  In the axisymmetric record problem, meridional Type-II records and critical signed compression are receiving branches rather than contradictions.  In the axisymmetric recurrence problem, first-use recharge has been forced toward typed noncompact boundary, age, frequency, source, and tail mechanisms, but the final changing-owner boundary-fed law has not been eliminated.

These endpoints are not versions of one theorem.  The general problem asks for singularity-preserving compactness or rigidity; the axisymmetric record problem asks for a receiving criterion; the axisymmetric recurrence problem asks for a same-solution first-use law.  A proof on one front may supply methodology for another, but no currency, norm, compactness theorem, or recurrence statement is identified across the two equations without an explicit bridge.

Accordingly, the companion papers prove structural and conditional theorems only in their stated scopes.  The present note adds no unconditional closure theorem.  It records the current logical architecture, including the places where compact methods stop and noncompact transition laws begin.  In particular, it does not prove unconditional regularity for three-dimensional Navier--Stokes, does not reduce the general problem to axisymmetric flow, and does not construct an axisymmetric singularity.

\end{document}